 \documentclass[10pt]{amsart}

\usepackage[all]{xypic}

\usepackage{latexsym}

\usepackage{amssymb}
\usepackage{amsfonts}
\usepackage{amscd}
\usepackage{amsmath,amsthm}
\usepackage{verbatim}

\newtheorem{lemma}{Lemma}[section]

\newtheorem{corollary}[lemma]{Corollary}

\newtheorem{lem}[lemma]{Lemma}
\newtheorem{prop}[lemma]{Proposition}
\newtheorem{thm}[lemma]{Theorem}
\newtheorem{cor}[lemma]{Corollary}

{

}

\theoremstyle{definition}

\theoremstyle{remark}

\numberwithin{equation}{section}

\newenvironment{pf}{\noindent{\bf Proof.}}{\hfill $\square$\medskip}

\def\PP{{\mathbb P}}

\def\ZZ{{\mathbb Z}}

\def\Fol{{\bar F}}

\def\Sol{{\bar S}}

\def\cTol{\overline{\cT}}

\def\0ol{{\bar 0}}
\def\1ol{{\bar 1}}
\def\2ol{{\bar 2}}
\def\ol2{{\bar 2}}
\def\3ol{{\bar 3}}
\def\4ol{{\bar 4}}
\def\5ol{{\bar 5}}
\def\6ol{{\bar 6}}
\def\7ol{{\bar 7}}
\def\8ol{{\bar 8}}
\def\9ol{{\bar 9}}

\def\bold0{{\bf 0}}
\def\bold1{{\bf 1}}
\def\bold2{{\bf 2}} 
\def\bold3{{\bf  3}}
\def\bold4{{\bf 4}}
\def\bold5{{\bf 5}}
\def\bold6{{\bf 6}}
\def\bold7{{\bf 7}}
\def\bold8{{\bf 8}}
\def\bold9{{\bf 9}}

\def\P2Skly{\PP^2_{Skly}}

\def\End{\operatorname {End}}
\def\Ext{\operatorname {Ext}}

\def\Gal{\operatorname {Gal}}

\def\Hom{\operatorname {Hom}}

\def\op{{\operatorname {op}}}

\def\Aut{\operatorname{Aut}}

\def\Div{\operatorname{Div}}

\def\End{\operatorname{End}}

\def\Ext{\operatorname{Ext}}

\def\gldim{\operatorname{gldim}}

\def\Gr{{\sf Gr}}

\def\Hom{\operatorname{Hom}}

\def\max{\operatorname{max}}

\def\mod{{\sf mod}}
\def\Mod{{\sf Mod}}

\def\op{\operatorname{op}}

\def\pdim{\operatorname{pdim}}

\def\Pic{\operatorname{Pic}}

\def\RHom{\operatorname {RHom}}
\def\RcHom{\operatorname {R\cHom}}

\def\Spec{\operatorname{Spec}}

\def\sup{\operatorname{sup}}

\def\ul1{\operatorname{\underline{1}}}

\def\a{\alpha}

\def\s{\sigma}

\def\sD{{\sf D}}

\def\wtM{{\widetilde{ M}}}

\def\cal{\mathcal}

\def\cE{{\cal E}}
\def\cF{{\cal F}}

\def\cI{{\cal I}}
\def\cJ{{\cal J}}

\def\cM{{\cal M}}

\def\cO{{\cal O}}

\def\cT{{\cal T}}

\def\cX{{\cal X}}

\def\coh{{\sf coh}}

\def\cHom{{\cal H}{\it om}}

\def\Qcoh{{\sf Qcoh}}

\def\dirlim{\mathop{\vtop{\baselineskip -100pt\lineskip -1pt\lineskiplimit 0pt
\setbox0\hbox{lim}\copy0\hbox to \wd0{\rightarrowfill}}}\limits}
\def\invlim{\mathop{\vtop{\baselineskip -100pt\lineskip -1pt\lineskiplimit 0pt
\setbox0\hbox{lim}\copy0\hbox to \wd0{\leftarrowfill}}}\limits}

\def\I11{{1 \kern -0.8pt \! \mbox{l}}}
\def\mumu{{\mu\kern-4.2pt\mu}}
\def\bfmu{{\mu\kern-4.2pt\mu}}
\def\2slash{\backslash \! \backslash}

\def\boxtimes{\setbox0\hbox{$\Box$}\copy0\kern-\wd0\hbox{$\times$}}

\newcommand{\mc}{\mathcal}
\newcommand{\struct}{\mc O}

\newcommand{\ang}[1]{\langle #1 \rangle}

\begin{document}

\title[A derived equivalence for a degree 6 del Pezzo surface]{A derived equivalence for a degree 6 del Pezzo surface over an arbitrary field}

\author{M. Blunk, S.J. Sierra, and S. Paul Smith}
\address{Department of Mathematics, Univ. of British Columbia.
\\ Department of Mathematics, Box 354350, Univ.  Washington, Seattle, WA 98195}
\email{mblunk@math.ubc.ca, sjsierra@math.washington.edu, smith@math.washington.edu}



\thanks{M. Blunk was supported by the National Science Foundation,
  Award No. 0902659.  S.J. Sierra was supported by the National
  Science Foundation, Award No. 0802935. S.P. Smith was supported by
  the National Science Foundation, Award No. 0602347}

\begin{abstract}
Let $S$ be a degree six del Pezzo surface over an arbitrary field
$F$. Motivated by the first author's classification of all such $S$ up to
isomorphism \cite{Bl} in terms of a separable $F$-algebra $B
\times Q \times F$, and by his K-theory isomorphism $K_n(S) \cong
K_n(B \times Q \times F)$ for $n \ge 0$, we prove an equivalence of
derived categories
$$
\sD^b(\coh S) \equiv \sD^b(\mod A)
$$ where $A$ is an explicitly given finite dimensional $F$-algebra
whose semisimple part is $B \times Q \times F$.
\end{abstract}

\maketitle
 
\section{Introduction}

We will work over an arbitrary field $F$. 

Throughout $S$ denotes a degree six del Pezzo surface over
$F$. Equivalently, $S$ is a smooth projective surface over $F$ whose
anti-canonical sheaf is ample and has self-intersection number 6.

Throughout $\Fol$ will denote a separable closure of $F$ and we will write
$$
\Sol=S_{\Fol}=S \times_{\Spec F}\Spec \Fol.
$$

In \cite{Bl}, the first author classified such $S$ up to isomorphism
by associating to $S$ a pair of separable $F$-algebras $B$ and $Q$, both defined
as endomorphism rings of certain locally free sheaves on
$S$. Furthermore, it was shown there that the algebraic K-theory of
$S$ is isomorphic to that of the algebra $B \times Q \times F$.

Let $\coh S$ denote the category of coherent sheaves on $S$ and let
$\mod A$ denote the category of noetherian right $A$-modules.  Let
$\equiv$ denote equivalence of derived categories.  Our main result
(Theorem~\ref{thm.main}) establishes a derived equivalence
\begin{equation}
\label{main.res}
\sD^b(\coh S) \equiv \sD^b(\mod A)
\end{equation}
where $A$ is a finite dimensional $F$-algebra whose semi-simple
quotient is $B \times Q \times F$.  We prove this equivalence by
constructing a tilting bundle $\cT$ on $S$ that has $A$ as its
endomorphism ring. (The definition of a tilting bundle is given in section 4.)
The main novelty of our approach is that we do not
make any assumptions on the base field $F$.  Since the field $F$ is
arbitrary, we cannot assume that $S$ is obtained by blowing up
$\PP^2_F$ (in fact $S$ could be a minimal surface), nor can we use
exceptional collections.
%
%
 
 \bigskip
 {\bf Acknowledgments.}  All three authors acknowledge the support of
 the National Science Foundation with gratitude. We would also like to
 thank the referee for their comments. This research was partially
 done while the first author was visiting the University of
 Washington, and he would like to thank the institution for its
 support and excellent working conditions. Finally, the first author
 would like to thank Aravind Asok, Baptise Calm\`{e}s, and Daniel
 Krashen for suggesting this problem.

\section{Basic facts about $\Sol$}

In this section, we give basic facts about the degree 6 del Pezzo
surface $\Sol$.  Since all the results here are well-known, we do not
give references.
    
There are six $(-1)$-curves on $\Sol$, which we may take to  lie in the following configuration:
 \begin{equation}
 \label{six.lines}
 \qquad
\UseComputerModernTips
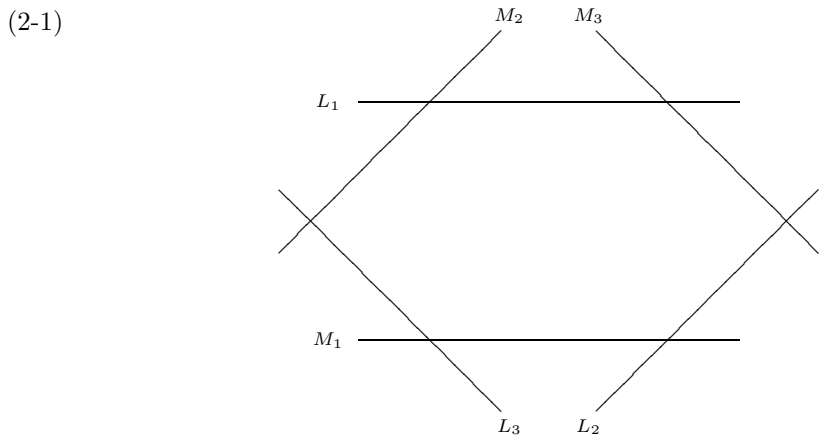
\xymatrix{
&&& &  \save []+<0cm,0.1cm>*\txt<4pc>{$\scriptstyle{M_2}$} \restore \ar@{-}[dddlll] 
& \ar@{-}[dddrrr]   \save []+<0cm,0.1cm>*\txt<4pc>{$\scriptstyle{M_3}$} \restore 
\\ 
&& \save []+<-0.3cm,0cm>*\txt<4pc>{$\scriptstyle{L_1}$} \restore 
\ar@{-}[rrrrr] &&&&&
\\  
&&&&&&&&
\\ 
&&&&&&&&
\\
&&
 \save []+<-0.3cm,0cm>*\txt<4pc>{$\scriptstyle{M_1}$} \restore 
 \ar@{-}[rrrrr] && &&&
\\
&&&& \ar@{-}[uuulll]   \save []+<0cm,-0.1cm>*\txt<4pc>{$\scriptstyle{L_3}$} \restore  
&  \save []+<0cm,-0.1cm>*\txt<4pc>{$\scriptstyle{L_2}$} \restore  \ar@{-}[uuurrr]  
}
\end{equation}
 The Picard group is
 $$ \Pic \Sol \cong {{\bigoplus_{i=1}^3 (\ZZ L_i \oplus \ZZ
     M_i)}\over{(M_i+L_j=M_j+L_i \; | \; 1\le i,j \le 3)}}.
 $$ Usually we only care about the class of a divisor in $\Pic
 \Sol$. We will write
 $$
 D_1  \sim D_2
 $$
 if $D_1$ and $D_2$ are linearly equivalent divisors.
 
 As remarked in the discussion after Prop. 2.1 in \cite{Bl}, the group
 of connected components of the group $\Aut \Sol$ is $S_2 \times S_3$,
 which can be identified with the automorphism group of the hexagon of
 $(-1)$-curves on $\Sol$. In particular, there is an element $\sigma
 \in \Aut(\Sol)$ that cyclically permutes the six exceptional
 lines. It is easy to see that $(1+\s)(1-\s^3)$ acts trivially on $\Pic
   \Sol$. 
 
An anti-canonical divisor is
$$
-K_{\Sol}:=L_1+L_2+L_3+M_1+M_2+M_3.
$$
This is ample. 
We define two particular divisors 
 \begin{equation}
 \label{H}
H\; := \;  L_1+M_2+M_3\; \sim\;  L_2+M_1+M_3 \; \sim \;  L_3+M_1+M_2
 \end{equation}
 and 
 \begin{equation}
 \label{H'}
H'\; := \;  L_1+L_2+M_3\; \sim \;  L_2+L_3+M_1 \; \sim \; L_3+L_1+M_2
 \end{equation}
 on $\Sol$.    Note that $\sigma(H) \sim H'$ and $\sigma^2(H) \sim H$.  

We define the {\sf degree} of a divisor $C$ on $\Sol$ as $\deg C = - C \cdot K$.
 Each exceptional line has degree 1.

There are two morphisms $f,f':\Sol \to \PP_\Fol^2$, each of which realizes
$\Sol$ as the blowup of $\PP_{\Fol}^2$ at three non-collinear
points. We choose these so that $f$ contracts the lines $L_1$, $L_2$,
and $L_3$ and $f'$ contracts the lines $M_1$, $M_2$, and $M_3$.  These
two morphisms induce injective group homomorphisms $f^*, f'^{*} : \Pic
\PP^2 \to \Pic \Sol$.  If $\ell$ is a line on $\PP^2_\Fol$, then $f^*\ell =
H$ and $f'^*\ell = H'$.

The action of $\Gal(\Fol/F)$ on the exceptional lines on $\Sol$
induces actions of $\Gal(\Fol/F)$ on
$$
\overline{\cI}:= \bigoplus_{i=0}^5  \cO_\Sol(\s^i H)
$$
and
$$
\overline{\cJ}:= \bigoplus_{i=0}^5  \cO_\Sol(\s^i( L_1+M_2))
$$ that are compatible with its action on $\Sol$.  In particular,
$\overline{\cI}$ and $\overline{\cJ}$ are $\Gal(\Fol/F)$-invariant.
It follows that the locally free sheaves $\overline{\cI}$ and
$\overline{\cJ}$ descend to locally free sheaves $\cI$ and $\cJ$ on
$S$.

 Define
 $$ \overline{ \cT}:= \overline{\cI} \oplus\overline{ \cJ} \oplus
 \cO_{\Sol}, \qquad \cT:= \cI \oplus \cJ \oplus \cO_S,
 $$
 and
 $$ B:= \End_{S} \cI, \qquad Q:= \End_{S} \cJ, \qquad
A:=\End_{S} \cT.
 $$ In \cite{Bl} it is shown that $S$ is determined up to isomorphism
by the pair of $F$-algebras $(B,Q)$.  (Actually, in \cite{Bl}, $B$ is
defined as $\big( \End_{S} \cI^\vee\big)^{\op}$.  Since sending a
homomorphism $\a:\cI \to \cI$ to its transpose $\a^\vee :\cI^\vee \to
\cI^\vee$ is an anti-isomorphism from $\End_{S} \cI$ to
$\End_{S} \cI^\vee$, our $B$ is the same as that in \cite{Bl},
and similarly for $Q$.)  As discussed in \cite{Bl}, the algebras $B$
and $Q$ are Azumaya over their centers, which are respectively \'{e}tale quadratic
and cubic extensions of $F$. Moreover, these \'{e}tale
centers can be recovered from the action of $\Gal(\Fol/F)$ on the
hexagon of $(-1)$-curves, as the action induces a 1-cocycle of
$\Gal(\Fol/F)$ with values in $S_2 \times S_3$, inducing a pair of
\'{e}tale extensions of $F$, quadratic and cubic.

We end this section with two results about the endomorphism algebra of $\cT$.

\begin{lem}
Let $A:=\End_{S} \cT$. Then
$$
A= \begin{pmatrix}
 B     &  \Hom_S(\cJ,\cI) & \Hom_S(\cO_S,\cI)   \\
   0   &  Q & \Hom_S(\cO_S,\cJ) \\
   0 & 0 & F 
\end{pmatrix}.
$$
\end{lem}
\begin{pf}
It suffices to show $\Hom_{\Sol}(\overline{\cI},\overline{\cJ}) =
\Hom_{\Sol}(\overline{\cI},\cO_\Sol) =
\Hom_{\Sol}(\overline{\cJ},\cO_\Sol) =0$.  However, each of these
three Hom-spaces is isomorphic to a direct sum of terms of the form
$H^0(\Sol, \cO_\Sol(D))$ for a divisor $D$ with $\deg D < 0$.  But if
$D$ has a section then $D \sim D'$ for some effective $D'$ so $\deg D
= -D'.K \geq 0$.  These Hom spaces are therefore zero.
 \end{pf}

The projective dimension of a left $T$-module is denoted by $\pdim_T
M$. The global homological dimension of $T$ is defined and denoted by
$$
\gldim T:= \sup\{\pdim_T M \; | \; M \in \Mod T\}.
$$

\begin{prop}
\label{prop.gldimA}
$\gldim A \le 2$.
\end{prop}
\begin{pf}
Let $R$ and $S$ be rings and $X$ an $R$-$S$-bimodule. 
 If $S$ is a semisimple ring, then 
$$
\gldim \begin{pmatrix} R & X \\ 0  & S \end{pmatrix} = \max \{ \pdim_R X +1, \gldim R \}.
$$
(See \cite[Prop. III.2.7]{ARS}.) Applying this result twice, first to
\begin{equation}
\label{temp.R}
A' := 
\begin{pmatrix} B &  \Hom(\cJ,\cI) \\ 0  & Q \end{pmatrix}
\end{equation}
then to $A$ with $R = A'$  and $S=F$, gives the desired result. 
\end{pf}

\section{Cohomology vanishing lemmas}
We will prove several results about vanishing of cohomology and
Ext-groups for sheaves on $S$.  These results will be used in
Section~\ref{sec-T} to show that $\cT$ is a tilting bundle and
therefore induces an equivalence of derived categories.

A key step in proving that $\cT$ is tilting is showing that
$\Ext^i_S(\cT,\cT)=0$ for $i>0$. This reduces, by flat base change, to
proving that $\Ext^i_\Sol(\overline{\cT},\overline{\cT})=0$. Given the
explicit description of $\overline{\cT}$ as a direct sum of invertible
sheaves, it suffices to prove that $h^1(D-D')=h^2(D-D')=0$ for all $D$
and $D'$ belonging to the list
\begin{equation}
\label{divisors}
H,\quad H',\quad  L_1+M_2,\quad L_2+M_3,\quad L_3+M_1, \quad 0.
\end{equation}
We will make repeated
 use of the relation $L_i + M_j \sim L_j + M_i$.

\begin{prop}
\label{prop.diffs} 
Let $D$ and $D'$ be divisors on $\Sol$ appearing in the list
(\ref{divisors}).  Then $$-3 \le \deg (D-D') \le 3.$$ Furthermore,
\begin{enumerate}
 \item
if $\deg (D-D')=1$, then $D-D'$ is linearly equivalent to an exceptional line. 
\item
if $\deg(D-D')=2$, then $D-D' \sim L_i+M_j$ for some $i \neq j \in
\{1,2,3\}$.
 \item
if $\deg(D-D')=3$, then $D-D'$ is linearly equivalent to either $H$ or
$H'$.
\item{} if $\deg(D-D')=0$, then $D-D'$ is linearly equivalent to
  either $0$, $L_i - L_j$, $L_i - M_i$, or $M_i-L_i$ for some $i,j \in
  \{1,2,3\}$.
\item{} if $\deg(D-D')<0$, then $D - D'$ is linearly equivalent to
  either $-L_i$, or $-M_j$, or $-L_i - M_j$, or $-H$, or $-H'$.
\end{enumerate}
\end{prop}
\begin{pf}
Exceptional lines have degree 1 so $\deg H=\deg H'=3$ and $\deg
(L_i+M_j)=2$.  It follows that the degree of $D-D'$ is between 3 and
-3.

(1) If $\deg(D-D')=1$, then $D$ is linearly equivalent to $H$ or $H'$
and $D' = L_i + M_j$ for some $i,j$. It follows from (\ref{H}) and
(\ref{H'}) that $D-D'$ is linearly equivalent to an exceptional line,
and every exceptional line can occur as $D-D'$.

(2) and (3) are obvious.

 (4)
In this case $D$ and $D'$ have the same degree. 

 If $\deg D=\deg D'=2$, then $D = L_i + M_j$ and $D' = L_k +
 M_\ell$. By considering all possible $i, j, k, \ell$, we see that $D - D'$
 is linearly equivalent to a divisor of the form $L_i - L_j$.

If $\deg D=\deg D'=3$, then, for example, $D \sim H$ and $D' \sim H'$,
and $D - D' \sim L_i -M_i$. Switching the roles of $H$ and $H'$, we see
$D - D' \sim M_i - L_i$.  Finally, we may have $D - D' \sim 0$.
  
 (5)
 This is the mirror of the cases (1)-(3).
\end{pf}

\begin{cor}
\label{cor.diffs}
Suppose $D$ is the difference of two divisors appearing in the list
(\ref{divisors}). If $\deg D \ge -2$, then there is an exceptional
line $E$ on $\Sol$ such that $D-E$ is also a difference of two
divisors appearing in the list (\ref{divisors}) and $D.E \ge -1$.
\end{cor}
\begin{pf}
This is established through case-by-case analysis using Proposition
\ref{prop.diffs} to look at all the possibilities for $D$.
\end{pf}

A divisor $D$ on $\Sol$ is {\sf good} if $h^1(D)=h^2(D)=0$.

  \begin{lemma}
  \label{lem.-H.good}
The divisors $-H$ and $-H'$ on $\Sol$ are good.
  \end{lemma}
  \begin{pf}
 The existence of the morphisms $f,f':\Sol \to \PP^2_\Fol$ allows us to use the
 Leray spectral sequence. The arguments for $-H$ and $-H'$ are the
 same so we only prove the result for $-H$.
 
Because $\Sol$ is a blowup of $\PP^2_\Fol$, $f_*\cO_{\Sol}=\cO_{\PP^2_{\Fol}}$ and
$R^jf_*\cO_{\Sol}=0$ if $j \ge 1$.
  
Since $\cO_{\Sol}(-H) \cong f^*\cO_{\PP^2_{\Fol}}(-\ell)$, the projection formula gives 
\begin{align*}
  R^jf_*\cO_{\Sol}(-H)  & =   R^jf_*\big( \cO_{\Sol} \otimes f^* \cO_{\PP^2_{\Fol}}(-\ell)\big)
  \\
  &
   \cong 
   R^jf_* \cO_{\Sol} \otimes  \cO_{\PP^2_{\Fol}}(-\ell) 
   \\
   & \cong 
   \begin{cases}
   \cO_{\PP^2_{\Fol}}(-\ell) & \text{if $j=0$}
   \\
   0 & \text{if $j \ne 0$}.
   \end{cases}
   \end{align*}
   The Leray spectral sequence
  $$
  H^i(\PP^2_\Fol,R^jf_*\cO_{\Sol}(-H)) \Rightarrow H^{i+j}(\Sol,\cO_{\Sol}(-H))
  $$
  therefore degenerates to give
  $$
  H^i(\Sol,\cO_{\Sol}(-H)) \cong H^i(\PP^2_\Fol,\cO_{\PP^2}(-\ell))
  $$
  for all $i$. The result follows because $H^i(\PP^2_\Fol,\cO_{\PP^2}(-\ell))=0$ for all $i$. 
  \end{pf}

\begin{lem}
\label{lem.good.induct}
Let $C$ be any divisor on $\Sol$, and let $E$ be one of the
$(-1)$-curves.  If $C-E$ is good and $C.E \ge -1$, then $C$ is good.
\end{lem}
\begin{pf}
The long exact sequence in cohomology associated to
$$
0 \to \cO_{\Sol}(C-E) \to \cO_{\Sol}(C) \to \cO_E(C) \to 0
$$
reads in part
$$
\UseComputerModernTips
\xymatrix{
\ar[r] &   H^1(\Sol,\cO_{\Sol}(C-E)) \ar[r] & H^1(\Sol,\cO_{\Sol}(C)) \ar[r] & H^1(\Sol, \cO_E(C)) \ar[r]& 
\\
\ar[r] & H^2(\Sol,\cO_{\Sol}(C-E)) \ar[r] & H^2(\Sol,\cO_{\Sol}(C))  \ar[r] & H^2(\Sol, \cO_E(C)).
}
 $$ By hypothesis, the left-most term in each row is zero. The
right-most term in each row is also zero because $H^i(\Sol,\cO_E(C))
\cong H^i(\PP^1_\Fol,\cO_{\PP^1}(C.E))$.
 Hence $C$ is good.
 \end{pf}

\section{The tilting bundle $\cT$}\label{sec-T}

In this section, we show that $\cT$ is a tilting bundle and prove our main result.

\begin{prop}
\label{prop.extTT}
Let $i \geq 1$. Then $\Ext^i_{S}( \cT,\cT) = 0$.
\end{prop}
\begin{pf}
By flat base change it suffices to prove this when $F$ is separably
closed so we assume that $F=\Fol$.  In that case $\Ext^i_S( \cT,\cT)$ is
isomorphic to a direct sum of terms of the form $H^i(S,\cO_S(D-D'))$
where $D$ and $D'$ are divisors in the list (\ref{divisors}).

It therefore suffices to show that $D-D'$ is good  
whenever $D$ and $D'$ are divisors in the list (\ref{divisors}).

We argue by induction on $\deg (D-D')$. By Proposition
\ref{prop.diffs}, $-3 \le \deg(D-D') \le 3$.  If $\deg (D-D')=-3$,
then $D-D'$ is good by Lemma \ref{lem.-H.good}. Now suppose that $ -2
\le \deg(D-D') \le 3$. By Corollary \ref{cor.diffs}, there is an
exceptional line $E$ such that $D-D'-E$ is a difference of divisors in
(\ref{divisors}) and $(D-D').E \ge -1$. By the induction
hypothesis, $D-D'-E$ is good, and it then follows from Lemma
\ref{lem.good.induct} that $D-D'$ is good.
\end{pf}

Since $\Sol$ is a del Pezzo surface of degree $\ge 6$ it is a toric variety
so we can, and will, make use of Cox's 
 homogeneous coordinate ring for it \cite{Cox}.

\begin{lem}
\label{lem.OY}
Every $\cF \in \coh \Sol$ has a
finite resolution in which all terms are direct sums of invertible
sheaves $\cO_{\Sol}(D)$ for various divisors $D$ on $\Sol$.
\end{lem}
\begin{pf}
Let $A$ be Cox's homogeneous coordinate ring for $\Sol$ \cite{Cox}. Then
$A$ is a polynomial ring with a grading by $\Pic (\Sol)$. Let $M$ be a
finitely generated graded $A$-module. Then $M$ has a finite projective
resolution in the category of graded $A$-modules. By \cite[Lemma
  2.2]{Sm}, every finitely generated projective graded $A$-module is a
direct sum of twists of $A$. The exact functor $\Gr(A,\Pic (\Sol)) \to
\Qcoh \Sol$, $M \rightsquigarrow \wtM$, described in
\cite[Thm. 3.11]{Cox} sends the resolution of $M$ to an exact sequence
in $\Qcoh \Sol$ in which the right-most term is $\wtM$ and all other
terms are direct sums of various $\cO_{\Sol}(D)$, $D \in \Div (\Sol)$. Given $\cF
\in \coh \Sol$, there is a finitely generated graded $A$-module $M$ such
that $\cF \cong \wtM$.
\end{pf} 

For the rest of this paper, we will work in the derived category.  If
$\sD$ is a triangulated category, we denote the shift of an object
$\cM$ by $\cM[1]$.  Recall that a subcategory of $\sD$ is {\em thick}
(\'epaisse) if it is closed under isomorphisms, shifts, taking cones of
morphisms, and taking direct summands of objects.

Let $\sD$ be a triangulated category and $\cE$ a set of objects in $\sD$. Then
\begin{itemize}
\item
$\sD^c$ denotes the full subcategory of $\sD$ consisting of the
  compact objects, i.e., those objects $C$ such that $\Hom_{\sD}(C,-)$
  commutes with direct sums;
  \item 
   $\langle \cE \rangle$ denotes the smallest thick full triangulated
subcategory of $\sD$ containing $\cE$;
  \item 
   $\cE^\perp$ denotes the full subcategory of $\sD$
consisting of objects $\cM$ such that $\Hom_{\sD}(E[i], \cM) = 0$ for
all $E \in \cE$ and all $i \in \ZZ$.
\end{itemize} 
We say that 
\begin{itemize}
  \item 
 $\cE$ {\sf generates} $\sD$ if $\cE^\perp = 0$ and that
 \item
 $\sD$ is {\sf compactly generated} if $(\sD^c)^\perp =0$. 
\end{itemize}
Clearly, if $\sD$ is compactly generated and $\langle \cE\rangle=\sD^c$, then $\cE$ generates $\sD$.

\begin{thm}
[Ravenel and Neeman \cite{Ne2}. Also see Thm. 2.1.2 in \cite{BvdB}]
\label{thm.RN}
Let $\sD$ be a compactly generated triangulated category.
Then a set of objects $\cE \subset \sD^c$ generates $\sD$ if and only if 
 $\langle \cE \rangle = \sD^c$.  
 \qed
\end{thm}

The unbounded derived categories $\sD(\Qcoh S)$ and $\sD(\Qcoh \Sol)$
are compactly generated. Moreover, $\sD(\Qcoh S)^c=\sD^b(\coh S)$ and
$\sD(\Qcoh \Sol)^c=\sD^b(\coh \Sol)$.
 
 \smallskip

{\bf Tilting bundles.}
Let $X$ be a projective scheme over a field $k$. A locally free sheaf $\cT \in \coh X$ is a
{\sf tilting bundle} if it generates $\sD(\Qcoh X)$ and $\Ext_{X}^i(\cT,\cT)=0$ for all $i>0$.

\begin{thm}
\label{thm.T.faithful}
$\cTol$ generates $\sD(\Qcoh \Sol)$ and $\langle \cTol \rangle =\sD^b(\coh \Sol)$. 
\end{thm}
\begin{pf}
By Theorem \ref{thm.RN}, it suffices to show that $\langle \cTol \rangle =\sD^b(\coh \Sol)$. 
Since $\langle \coh \Sol \rangle= \sD^b(\coh \Sol)$ it suffices to
show that every coherent $\cO_{\Sol}$-module belongs to $\langle \cTol   \rangle$.
 
  If $D$ is an effective divisor on $\Sol$ we write $\cI_D$ for the ideal
  vanishing on $D$ as a scheme. Thus $\cI_D \cong \cO_{\Sol}(-D)$. Whenever
  we write an arrow $\cO_{\Sol}(-D) \to \cO_{\Sol}$ it will be with the tacit
  understanding that this is the composition of an isomorphism
  $\cO_{\Sol}(-D) \to \cI_D$ followed by the inclusion $\cI_D \to \cO_{\Sol}$.
  
 Since $M_3 \cdot(L_1+M_2+M_3)=0$,   $\cO_{M_3} \cong \cO_{M_3}(L_1+M_2+M_3)$. It follows from the 
 exact sequences
 $$
 0 \to \cO_{\Sol}(L_1+M_2) \to \cO_{\Sol}(L_1+M_2+M_3) \to \cO_{M_3}(L_1+M_2+M_3) \to 0
 $$
 and
 $$
 0 \to \cO_{\Sol}(-M_3) \to \cO_{\Sol} \to \cO_{M_3} \to 0
 $$
 that $\cO_{M_3}$ and $\cO_{\Sol}(-M_3)$ belong to $\langle \cTol \rangle$.
 Hence  $\cO_E$ and $\cO_{\Sol}(-E)$ belong to $\langle \cTol \rangle$ for all exceptional lines $E$.

Since $L_i .L_k = 0$ if $i \neq k$, there is an exact sequence
\[ 0 \to \struct_{\Sol}(-L_i - L_k) \to \struct_{\Sol}(-L_i) \oplus \struct_{\Sol}(-L_k) \to \struct_{\Sol} \to 0.\]
Twisting by $L_i + M_j +L_k$, we obtain
\[ 0 \to \struct_{\Sol}(M_j) \to \struct_{\Sol}(M_j +L_k) \oplus \struct_{\Sol}(L_i + M_j) \to \struct_{\Sol}(L_i +M_j+L_k) \to 0.\]
Therefore, $\struct_{\Sol}(M_j) \in \ang \cTol$.
From the exact sequence
\[ 0 \to \struct_{\Sol} \to \struct_{\Sol}(M_j) \to \struct_{M_j}(M_j) \to 0,\]
we deduce that $\struct_{M_j}(M_j) \in \ang{\cTol}$.  

It follows that $\cO_E(E) \in  \langle \cTol \rangle$ for every exceptional curve $E$. But $\cO_E$ is also in 
$  \langle \cTol \rangle$ so,  because $\sD^b(\coh \PP^1_{\Fol})$ is generated by 
  $ \cO_{\PP^1_{\Fol}}$ and $ \cO_{\PP^1_{\Fol}}(-1)$,  it follows that  $\sD^b(\coh E) \subset  \langle \cTol \rangle$.
  Hence $\cO_E(D) \in  \langle \cTol \rangle$ for all divisors $D$ on $\Sol$. 
 
Suppose $\cO_{\Sol}(D) \in   \langle \cTol \rangle$. Then $\cO_{\Sol}(D-E) \in   \langle \cTol \rangle$ because there is an 
 exact sequence
$$
0 \to \cO_{\Sol}(D-E) \to \cO_{\Sol}(D) \to \cO_E(D) \to 0.
$$
Likewise, $\cO_{\Sol}(D+E) \in   \langle \cTol \rangle$ because  there is an exact sequence
$$
0 \to \cO_{\Sol}(D) \to \cO_{\Sol}(D+E) \to \cO_E(D+E) \to 0.
$$ It follows that $ \langle \cTol \rangle$ contains $\cO_{\Sol}(D)$
for all $D \in \Div \Sol$ and therefore, by Lemma \ref{lem.OY},
contains $\cF$ for every $\cF \in \coh \Sol$.
\end{pf}

When $F$ is not separably closed $\cT$ need not split as a direct sum of line bundles so the arguments
in Theorem \ref{thm.T.faithful} can not be used to prove directly that $\langle \cT \rangle = \sD^b(\coh S)$. 
Instead we will show that $\cT$ generates $\sD(\Qcoh S)$ and then apply Theorem~\ref{thm.RN}.  

\begin{thm}\label{thm.main}
Let $F$ be an arbitrary field. Then
 $$\RHom_S(\cT,-):\sD^b(\coh S) \to \sD^b(\mod A)$$ is an equivalence of categories.
 \end{thm}
 \begin{pf}
We will show that $\cT$ generates $\sD(\Qcoh S)$. It will then follow from  Theorem \ref{thm.RN} that
$$
\langle \cT \rangle = \sD(\Qcoh S)^c =\sD^b(\coh S).
$$
By Proposition~\ref{prop.extTT}, $\Ext^i_S(\cT, \cT) = 0$ for $i>0$.  By Proposition~\ref{prop.gldimA}, $A = \End_S(\cT)$ has finite global dimension. 
Thus we have shown that $\cT$ is a tilting bundle and our theorem will then follow directly from 
\cite[Thm. 3.1.2]{Baer} (or \cite[Thm. 7.6]{HVdB}). 

Let $\cM \in \sD(\Qcoh S)$  and suppose $\RHom_S(\cT,\cM)=0$. We must show that $\cM =0$.
 
 Since $\cT$ is locally free, $\cHom_{S}(\cT,-)$ and $\cT^\vee \otimes_{S}-$ are exact 
 functors on $\Qcoh S$. Likewise, $\cHom_{\Sol}(\cTol,-)$ and $\cTol^\vee \otimes_{\Sol}-$ are exact 
 functors on $\Qcoh \Sol$. Thus, for example, $\RcHom_{S}(\cT,\cM)$ 
 can be computed on $\sD(\Qcoh S)$
 by applying $\cHom_{S}(\cT,-)$ to each individual term in $\cM$.

 Consider the cartesian square
 $$
 \UseComputerModernTips
\xymatrix{
\Sol \ar[rr]^v \ar[d]_q && S \ar[d]^p
\\
\Spec (\Fol) \ar[rr]_u && \Spec(F).
}
$$
Since $u$ (and therefore $v$) is flat, 
 the natural transformation
$$
u^* Rp_*\to Rq_*v^*
$$
is an isomorphism of functors from $\sD(\Qcoh S)$ to $\sD( \Fol)$ \cite[(3.18)]{Huy}. We now have 
\begin{align*}
0 \; = \; u^* \RHom_S(\cT,\cM) &  \cong \; u^*Rp_* \RcHom_S(\cT,\cM) 
\quad \hbox{by \cite[p.85]{Huy}}
\\
&  \cong \; Rq_*v^* \RcHom_S(\cT,\cM)  
\quad \hbox{by \cite[(3.18)]{Huy}}
\\
&  \cong \; Rq_*v^*  (\cT^\vee \otimes^L_{S} \cM ) 
\\
&  \cong \; Rq_*(\cTol^\vee \otimes^L_{\Sol} Lv^* \cM ) 
\\
&  \cong \; Rq_*\RcHom_{\Sol}(\cTol,  Lv^* \cM ) 
\\
&  \cong \;  \RHom_{\Sol}(\cTol,  Lv^* \cM ) .
\end{align*}
 But $\cTol$ generates $\sD(\Qcoh \Sol)$ so $v^*\cM = 0$. Since
$v^*$ is faithful, $\cM=0$, and we are done.
\end{pf}

  \begin{corollary}[cf. \cite{Bl}, Corollary 5.2]
The functor $\Hom_S (\cT, -): \coh (S) \rightarrow \mod A$ induces an
isomorphism 
\[\Hom_S(\cT, -): K_*(S) \rightarrow K_*(F \times B \times Q).\]
  \end{corollary}
\begin{pf}
It follows from Theorem $1.98$ of \cite{ThTr} that the equivalence of
derived categories found in Theorem \ref{thm.main} induces an
isomorphism in $K$-theory
\[ \Hom_S(\cT, -): K_* (\coh S) \rightarrow K_* (\mod A).\]
Moreover, $A$ has a nilpotent ideal $I$ so that $A/I$ is isomorphic to
its semi-simple quotient $F \times B \times Q$. Thus, it follows that
the $K$-theory of $A$ is isomorphic to that of $F \times B \times Q$,
and we recover the isomorphism found in \cite{Bl}.
\end{pf}

\end{document}